\documentclass[12pt,letterpaper]{article}
\usepackage{amsfonts}
\usepackage{amssymb}
\usepackage{latexsym}  

\let\SavedRightarrow=\Rightarrow
\usepackage{marvosym}
\let\Rightarrow=\SavedRightarrow

\usepackage{amsmath}   

\DeclareMathAlphabet{\mathbbold}{U}{bbold}{m}{n}

\renewcommand\AA{{\mathcal A}}
\newcommand\PP{{\mathcal P}}
\newcommand\PS{\mathcal{PS}}

\newcommand\RF{\mathsf{REQ}}

\newcommand\RRR{{\mathbb R}}
\newcommand\NNN{{\mathbb N}}
\newcommand\TTT{{\mathbb T}}
\newcommand\CCC{{\mathbb C}}
\newcommand\PPP{{\mathbb P}}

\newcommand\DDDD{{\mathfrak D}}
\newcommand\PPPP{{\mathfrak P}}

\newcommand\cchi{{\raise 2 pt \hbox{$\chi$}}}

\newcommand\one{\mathbbold{1}} 
\newcommand\cat{^{\mathord{\frown}}}  

\newcommand\res{\mathord {\upharpoonright}}  

\newcommand\comp{\mathrm{comp}}  
\newcommand\dom{\mathrm{dom}}  
\newcommand\diam{\mathrm{diam}}   
\newcommand\cdisc{\overline{D}}  

\newcommand\lea{\sqsubseteq}  
\newcommand\leac{\, \mbox{\rlap{$\sqsubseteq$}  %
   \raise 2.3pt \hbox{\scriptsize $\mathrm{c}$}}\; }

\newcommand\iv{^{-1}} 

\newcommand\onto{\twoheadrightarrow}

\newcommand\eop{\raisebox{-2pt}{\mbox{\Huge \Smiley}}}

{\end{enumerate}}
\newenvironment{itemizz}{\begin{itemize}\setlength{\itemsep}{-1mm}} %
{\end{itemize}}                              

\newtheorem{theorem}{Theorem}[section]
\newtheorem{definition}[theorem]{Definition}
\newtheorem{lemma}[theorem]{Lemma}
\newtheorem{corollary}[theorem]{Corollary}

\newenvironment{proof}{{\bf Proof.}}{\eop\medskip}
\newenvironment{proofof}[1]{\medskip \textbf{Proof of #1.}}{\eop\medskip}

\begin{document}

\title{Inverse Limits and Function Algebras}

\author{Joan E. Hart\footnote{University of Wisconsin, Oshkosh,
WI 54901, U.S.A.,
\ \ hartj@uwosh.edu}
\  and
Kenneth Kunen\footnote{University of Wisconsin,  Madison, WI  53706, U.S.A.,
\ \ kunen@math.wisc.edu}
\thanks{Both authors partially supported by NSF Grant DMS-0097881.}
}

\maketitle

\begin{abstract}
Assuming
Jensen's principle $\diamondsuit$, there is a compact Hausdorff space $X$
which is hereditarily Lindel\"of, hereditarily separable, and connected,
such that no closed subspace of $X$ is both perfect and totally disconnected.
The Proper Forcing Axiom implies that there is no such space.
The $\diamondsuit$ example also fails to satisfy the
CSWP (the complex version of the Stone-Weier\-strass Theorem).
This space cannot contain the two earlier examples of
failure of the CSWP, which were totally disconnected ---
specifically, the Cantor set (W. Rudin) and $\beta\NNN$ (Hoffman and Singer).
\end{abstract}

\section{Introduction} 
\label{sec-intro}
All topologies discussed in this paper are assumed to be Hausdorff.
It is well-known that if $X$ is compact and second countable
and not scattered, then $X$ has a subspace homeomorphic
to the usual Cantor set, $2^\omega$.
This is not true of non--second countable spaces.
For example, the double arrow space of Alexandroff and Urysohn \cite{AU}
is  compact and not scattered, but is only first countable and 
does not contain a Cantor subset.

The double arrow space is 
also HS (hereditarily separable) and HL (hereditarily Lindel\"of);
that is, all subspaces are both separable and Lindel\"of
(see \cite{ENG} Exercise 3.10.C).
It is also a LOTS; that is, a totally ordered set with
its order topology.
The double arrow space is also totally disconnected, and it is
natural to ask whether there is a connected version of it.
This turns out to be independent of ZFC.  Under the Proper
Forcing Axiom (PFA), there is no such space:

\begin{theorem}
\label{thm-PFA-HL}
Assuming PFA, every compact HL space is either totally disconnected
or contains a copy of the Cantor set.
\end{theorem}

On the other hand, by Theorem \ref{thm-weird-HL},
there will be such a space assuming
Jensen's principle $\diamondsuit$, which is true in
G\"odel's universe of constructible sets.

\begin{definition}
\label{def-weird}
A space $X$ is \emph{weird} iff $X$ is compact and not scattered,
and there is no $P \subseteq X$ such that $P$ is perfect
and totally disconnected.
\end{definition}

So, $X$ cannot be second  countable.  However,

\begin{theorem}
\label{thm-weird-HL}
Assuming $\diamondsuit$, there is a weird $X$ such that
$X$ is HS and HL.
\end{theorem}

Note that a compact $X$ is HL iff every closed set is a $G_\delta$
(see \cite{ENG} Exercise 3.8.A(c)).
Applying this to the points, we see that $X$ must 
be first  countable.

We can also get our $X$ to fail the complex version of
the \emph{Stone-Weierstrass Theorem}.
This theorem involves subalgebras
of $C(X,\RRR)$, and is true for all compact $X$.
If one replaces the real numbers $\RRR$ by the complex numbers $\CCC$,
the ``theorem'' is true for some $X$ and false for others,
so it becomes a \emph{property} of $X$:

\begin{definition}
\label{def-cswp}
If $X$ is compact, then
$C(X) = C(X,\CCC)$ is the algebra of continuous
complex-valued functions on $X$,
with the usual supremum norm.
$\AA \lea C(X)$ means that $\AA$ is a subalgebra of $C(X)$
which separates points and contains the constant functions.
$\AA \leac C(X)$ means that $\AA \lea C(X)$ and
$\AA$ is closed in $C(X)$.
$X$ has the \emph{Complex Stone-Weierstrass Property (CSWP)}
iff every $\AA \lea C(X)$ is dense in $C(X)$.
\end{definition}

Classical examples from the 1800s show that the CSWP is false
for many $X$.  In particular:

\begin{definition}
\label{def-disc}
$D$ denotes the open unit disc in $\CCC$ and $\TTT$ denotes the unit circle.
The \emph{disc algebra} $\DDDD \leac C(\cdisc)$ is the set of $f \in C(\cdisc)$
which are holomorphic on $D$.
\end{definition}

Then $\DDDD$ refutes the CSWP of $\cdisc$, and
$\DDDD\res\TTT = \{f \res \TTT : f \in \DDDD\}$ refutes
the CSWP of $\TTT$.
Further negative results were obtained in 1956
by Rudin \cite{RUD1} and in 1960 by Hoffman and Singer \cite{HS}
(see also \cite{HO}):
\begin{itemizz}
\item[1.] \cite{RUD1}  Every compact $X$ containing a copy of the
Cantor set fails the CSWP.
\item[2.] \cite{HS}  Every compact $X$ containing a copy of $\beta\NNN$
fails the CSWP.
\end{itemizz}
Actually, \cite{HS} does not mention $\beta\NNN$, and used instead
$S =$ the Stone space of a separable measure algebra,
but that is equivalent, since each of $S$ and $\beta\NNN$ contains 
a copy of the other.
The first non-trivial positive result is due to Rudin \cite{RUD2}, 
and some more recent positive results are contained in 
\cite{HK,KU}.  In particular,
\begin{itemizz}
\item[3.] \cite{RUD2}  Every compact scattered space
satisfies the CSWP.
\item[4.] \cite{KU}  Every compact LOTS
which does not contain a copy of the Cantor satisfies the CSWP.
\end{itemizz}
By (4), the double arrow space 
is an example of a non-scattered space which has the CSWP.
Results (1) through (4) might suggest the (highly unlikely) conjecture that
a compact $X$ has the CSWP whenever it contains neither $\beta\NNN$
nor a Cantor set. Under $\diamondsuit$, this is refuted by:

\begin{theorem}
\label{thm-main}
Assuming $\diamondsuit$, there is a weird $X$ such that
$X$ is HS and HL and $X$ fails the CSWP.
\end{theorem}

As Rudin pointed out, (1)(3) imply that for $X$ compact metric,
$X$ has the CSWP iff $X$ does not contain a Cantor subset.
By (1)(4), the same ``iff'' holds when $X$ is a compact LOTS\@.
By (2), the ``iff'' does not hold for arbitrary compact spaces,
but one might hope to prove it for
some other spaces which are small in some way.
Theorem \ref{thm-main} puts some bounds on this hope.

Obviously, Theorem \ref{thm-main} implies Theorem \ref{thm-weird-HL},
but we shall prove Theorem \ref{thm-weird-HL} first.
We then explain what needs to be added to the construction
to obtain Theorem \ref{thm-main}.
Both proofs are essentially inverse limit constructions.
For Theorem \ref{thm-weird-HL}, we obtain $X \subset [0,1]^{\omega_1}$
by an inductive construction; at stage $\alpha < \omega_1$,
we determine the projection, $X_\alpha$, of $X$ on $[0,1]^\alpha$.
Then, $X$ may be viewed as the inverse limit of
$\langle X_\alpha : \alpha < \omega_1 \rangle$.
For Theorem \ref{thm-main}, we replace $[0,1]$ by $\cdisc$.

Theorem \ref{thm-weird-HL} is proved in
Section \ref{sec-weird}, which also gives some more information
about weird spaces.
Theorem \ref{thm-main} is proved in Section \ref{sec-poly},
using a fact about peak points proved in Section \ref{sec-peak}.
Theorem \ref{thm-PFA-HL} is proved in Section \ref{sec-forcing},
which may be read immediately after Section \ref{sec-weird}.

\section{Weird Spaces}
\label{sec-weird}

We list some easy properties of weird spaces:

\begin{definition}
$\comp(x,X)$ denotes the connected component of the point $x$
in the space $X$.
\end{definition}

\begin{lemma}
\label{lemma-weird-easy}
If $X$ is weird then:
\begin{itemizz}
\item[1.] If $Y \subseteq X$ and $Y$ is closed, then $Y$ is either
scattered or weird.
\item[2.] For some $x \in X$:  $\comp(x,X)$ is not a singleton,
so that $\comp(x,X)$ is weird and connected.
\item[3.] $X$ is not second countable.
\item[4.] $X$ is not a LOTS.
\end{itemizz}
\end{lemma}
\begin{proof}
For (4), if $X$ is a LOTS, let $S\subset X$ be countable and
order-isomorphic to the rationals.  Since $\overline S$ cannot be
totally disconnected, it contains an interval isomorphic
to the unit interval in $\RRR$, contradicting (1) and (3).
\end{proof}

We shall see that no subspace of a countable product of LOTS can
be weird either.  First:

\begin{lemma}
\label{lemma-weird-onto}
If $X$ is weird and $f$ maps $X$ continuously onto $Y,$
then either $Y$ is weird or some $f\iv\{y\}$ is weird.
\end{lemma}
\begin{proof}
Assume no $f\iv\{y\}$ is weird.  Then
each $f\iv\{y\}$ is scattered.

Note that $Y$ cannot be scattered.  To see this, let $K$ be the perfect
kernel of $X$.  If $y$ is an isolated point of $f(K)$, then
$K \cap f\iv\{y\}$ is scattered and clopen in $K$, a contradiction.

If $Y$ is not weird, fix 
$P \subseteq Y$ such that $P$ is  perfect
and totally disconnected.
Then for $x \in f\iv(P)$, $\comp(x, f\iv(P)) \subseteq f\iv \{f(x)\}$,
which is scattered, so $\comp(x, f\iv(P)) = \{x\}$.
Thus, $f\iv(P)$ is totally disconnected, and hence scattered
(since  $X$ is weird), which is a contradiction, since
$P = f(f\iv(P))$ is not scattered.
\end{proof}

\begin{corollary}
\label{cor-weird-finite}
Suppose that $X$ is weird and
$X \subseteq \prod_{j<n} Z_j$, where $n$ is finite and
each $Z_j$ is compact.  Then some $Z_j$ has a weird subspace.
\end{corollary}
\begin{proof}
Induct on $n$, using Lemma \ref{lemma-weird-onto}.
\end{proof}

We now prove the same result for countable products.
First, we introduce some notation for products and projections:

\begin{definition}
\label{def-proj}
If $Z_\xi$ are spaces for $\xi  < \beta$ 
then
$\pi^\beta_\alpha : \prod_{\xi < \beta} Z_\xi \onto  \prod_{\xi < \alpha}Z_\xi$
\textup(for $\alpha \le \beta$\textup) and
then
$\varphi^\beta_\alpha : \prod_{\xi < \beta} Z_\xi \onto  Z_\alpha$
\textup(for $\alpha <\beta$\textup) and
are the natural projections.
If
$\vec z = \langle z_\xi : \xi < \beta \rangle \in  \prod_{\xi < \beta} Z_\xi $,
then $\varphi^\beta_\alpha(\vec z) = z_\alpha$ and
$\pi^\beta_\alpha(\vec z) = \langle z_\xi : \xi < \alpha \rangle$.
We sometimes write $\varphi_\alpha$
for $\varphi^\beta_\alpha$ when $\beta$ is clear from context.
\end{definition}

\begin{lemma}
Suppose that $X$ is weird and
$X \subseteq \prod_{j<\omega} Z_j$, where 
each $Z_j$ is compact.  Then some $Z_j$ has a weird subspace.
\end{lemma}
\begin{proof}
Assume that no $Z_j$ has a weird subspace; we shall derive
a contradiction.

Let $X_n = \pi^\omega_n(X) \subseteq  \prod_{j<n} Z_j$.
By Corollary \ref{cor-weird-finite}, no $X_n$ has a weird subspace.

View $\bigcup_n X_n$ as a tree, where $X_n$ is the 
$n^\mathrm{th}$ level, and the tree order $<$ satisfies 
$y < z$ iff $y = \pi^n_m(z)$ whenever $m < n$,
$y \in X_m$ and $z \in X_n$.
Let $W_n$ be the set of all $y \in X_n$ such that
$X \cap (\pi^\omega_n)\iv\{y\}$ is weird (equivalently, non-scattered).
Note that $\bigcup_n W_n$ is a subtree of $\bigcup_n X_n$;
equivalently, $\pi^n_m(W_n) \subseteq W_m$ whenever $m < n$.

First, note that if $P \subseteq X$ is closed and not scattered,
then $W_n \cap \pi^\omega_n(P) \ne \emptyset$ for each $n$.
To see this, use the fact that $P$ is weird and
$\pi^\omega_n(P)$ is not weird, and
apply Lemma \ref{lemma-weird-onto} to $\pi^\omega_n \res P$.

It follows that 
$\bigcup_n W_n$ is a perfect tree; that is,
if $y \in W_m$, then for some $n > m$, there are more than one
$z \in W_n$ such that $\pi^n_m(z) = y$.
To see this, let $P_0,P_1$ be disjoint perfect subsets of 
$X \cap (\pi^\omega_m)\iv\{y\}$, and choose $n$ such that 
$\pi^\omega_n(P_0) \cap \pi^\omega_n(P_1) = \emptyset$.
If $z_\ell \in W_n \cap \pi^\omega_n(P_\ell)$ (for $\ell = 0,1$),
then $z_\ell \in W_n$ and $\pi^n_m(z_\ell) = y$ and
$z_0 \ne z_1$.

But now we can choose a Cantor subtree.  That is, we can choose
finite non-empty $F_n \subseteq W_n$ so that 
$m < n \to \pi^n_m(F_n) = F_m$ and
for each $m$, there is an $n > m$ such that 
$|F_n \cap (\pi^n_m)\iv\{y\} | \ge 2$ for all $y \in F_m$.
Then $\{x \in X : \forall n [ \pi^\omega_n(x) \in F_n ] \}$ is
homeomorphic to the Cantor set, a contradiction.
\end{proof}

In particular, by Lemma \ref{lemma-weird-easy}, and the observation
that every closed subspace of a compact LOTS is a compact LOTS:

\begin{corollary}
\label{cor-weird-countable}
Suppose that $X \subseteq \prod_{j<\omega} Z_j$, where 
each $Z_j$ is compact and is either second countable or a LOTS.
Then $X$ is not weird.
\end{corollary}

We now turn to a proof of Theorem \ref{thm-weird-HL},
which obtains a weird subspace  of
an uncountable product, $[0,1]^{\omega_1}$.
There are many such constructions in the literature;
we follow the specific approach in \cite{DK}\S4,
which uses
irreducible projections (see \cite{ENG} Exercise 3.1.C)
to ensure that the space is HS and HL.

\subsection{The Construction}
\label{sec-construct}

We shall get $X  = X_{\omega_1} \subseteq [0,1]^{\omega_1}$ with
$X_\alpha = \pi^{\omega_1}_\alpha(X) \subseteq [0,1]^{\alpha}$ 
satisfying:
\begin{itemizz}
\item[0.] $X_1 = [0,1]$.
\item[1.] $X_\alpha$ is connected whenever $1 \le \alpha \le \omega_1$.
\item[2.] $\pi^\beta_\alpha: X_\beta \onto X_\alpha$ is
irreducible whenever $1 \le \alpha \le \beta \le \omega_1$.
\end{itemizz}
In particular, $\pi^{\omega_1}_1: X_{\omega_1} \onto X_1$ will be
irreducible, so $X$ will be separable and have no isolated points.  To make $X$ 
HS, we get $P_\alpha$ and $\PP_\alpha$ for $1 \le \alpha < \omega_1$ so that:
\begin{itemizz}
\item[3.] $\PP_\alpha$ is a countable family of closed subsets 
of $X_\alpha$ and $P_\alpha \in \PP_\alpha$.
\item[4.] If $P \in \PP_\alpha$, then 
\begin{itemizz}
\item[a.]
$\pi^{\alpha+1}_\alpha : (\pi^{\alpha+1}_\alpha)\iv(P) \onto P$
is irreducible, and
\item[b.]
$(\pi^{\beta}_\alpha)\iv(P) \in \PP_{\beta}$ whenever
$\alpha < \beta < \omega_1$.
\end{itemizz}
\end{itemizz}

\begin{lemma}
Requirement $(4)$ implies  that 
$\pi^\beta_\alpha : (\pi^\beta_\alpha)\iv(P) \onto P$ is irreducible
for all $\beta \ge \alpha$ whenever $\alpha \le \beta \le \omega_1$.
\end{lemma}
\begin{proof}
Induct on $\beta$.
\end{proof}

To get $X$ to be HL and HS, we add the next requirement:

\begin{itemizz}
\item[5.] If $F \subseteq X$ is closed, then
$\pi^{\omega_1}_\alpha(F) = P_\alpha$ for some $\alpha < \omega_1$.
\end{itemizz}

\begin{lemma}
\label{lem-HLHS}
Requirements $(4)(5)$ imply that $X$ is HL and HS.
\end{lemma}
\begin{proof}
To see that $X$ is HL,  use (5) and (4) to see that
every closed $F \subseteq X$ is a $G_\delta$: 
For every closed subset $F$ of $X$,  we have 
$\pi^{\omega_1}_\alpha(F) = P \in \PP_\alpha$.
Then by irreducibility, $F = (\pi^{\omega_1}_\alpha)\iv (P)$,
so that $F$ is a $G_\delta$.
Also by (5), all closed $F \subseteq X$ are separable,
so $X$ is HS (since it is HL and hence first countable).
\end{proof}

Conditions (0)--(5) are consistent with all $\pi^\beta_\alpha$
being homeomorphisms, which would make $X$ homeomorphic to $[0,1]$.
To make $X$ weird, we also choose $h_\alpha, p_\alpha$, and $q_\alpha^n$
for $n < \omega$ and $0 < \alpha < \omega_1$ so that:

\begin{itemizz}
\item[6.] $p_\alpha \in X_\alpha$ and
$h_\alpha \in C(X_\alpha \backslash \{p_\alpha\},\; [0,1])$
and $X_{\alpha+1} = \overline{h_\alpha}$.
\item[7.] $q_\alpha^n \in  X_\alpha \backslash\{p_\alpha\}$, and
$\langle q_\alpha^n: n \in \omega \rangle \to p_\alpha$, and all points
of $[0,1]$ are limit points of
$\langle h_\alpha(q_\alpha^n): n \in \omega \rangle$, and
$\{p_\alpha\} \times [0,1] \in \PP_{\alpha+1}$.
\item[8.] For each $P \in \PP_\alpha$, either $p_\alpha \notin P$, 
or $p_\alpha \in P$ and $q_\alpha^n \in P$ for all but finitely many $n$.
\end{itemizz}
As usual, we identify $h_\alpha$ with its graph, which is a subset
of $X_\alpha \times [0,1]$;
we also identify $[0,1]^\alpha \times [0,1]$ with $[0,1]^{\alpha+1}$.

\begin{lemma}
Requirements $(0)(6) (7)$ imply requirements $(1) (2)$.
\end{lemma}
\begin{proof}
Induct on $\alpha$.
By (6), $\pi^{\alpha+1}_\alpha : X_{\alpha+1} \onto X_\alpha$
is one-to-one at all points not in $(\pi^{\alpha+1}_\alpha)\iv\{p_\alpha\}$.
The first part of (7) implies that
$\{p_\alpha\} \times [0,1] \subseteq X_{\alpha+1}$.
\end{proof}

\begin{lemma}
\label{lemma-nine}
Requirements $(0)$ -- $(8)$ imply that if $P \in \PP_\alpha$ is connected, then
$(\pi^{\omega_1}_\alpha)\iv(P)$ is connected.
\end{lemma}
\begin{proof}
By (8), if $P \in \PP_\alpha$ is connected, then 
$(\pi^{\alpha+1}_\alpha)\iv(P)$ will also be connected, so 
now prove that $(\pi^{\beta}_\alpha)\iv(P)$ is connected
by induction on $\beta \le \omega_1$.
\end{proof}

To help make $X$ weird we add the requirement:

\begin{itemizz}
\item[9.] If $F \subseteq X$ is closed and not scattered,
then for some $\alpha < \omega_1$,
$\pi^{\omega_1}_\alpha(F) = P_\alpha$
and $P_\alpha$ is not scattered and $p_\alpha \in P_\alpha$.
\end{itemizz}

Note that we cannot simply omit (5) in favor of (9),
since Lemma \ref{lem-HLHS} uses (5) for all closed $F$,
including singletons.

\begin{lemma}
Requirements $(0)$ -- $(9)$  imply that $X$ is weird.
\end{lemma}
\begin{proof}
By (9), every closed non-scattered $F\subseteq X$ satisfies
$\pi^{\omega_1}_\alpha(F) = P_\alpha$, 
for some $\alpha < \omega_1$,
with $P_\alpha$ not scattered and $p_\alpha \in P_\alpha$.
Such $F$ therefore contain 
$(\pi^{\omega_1}_{\alpha+1})\iv(\{p_\alpha\} \times [0,1])$.
By  (7) and Lemma \ref{lemma-nine}, each 
$(\pi^{\omega_1}_{\alpha+1})\iv(\{p_\alpha\} \times [0,1])$
is a connected subspace of $X$.
\end{proof}

\begin{proofof}{Theorem \ref{thm-weird-HL}}
To get (5) and (9), use $\diamondsuit$ to capture all
closed subsets of $[0,1]^{\omega_1}$.  
To get (7)(8) for a fixed $\alpha$:  First, list $\PP_\alpha$
as $\{Q^n : n \in \omega\}$, with $Q^0 = P_\alpha$.
Let $d$ be a metric on $X_\alpha$.
Choose perfect $F^n \subseteq X_\alpha$ for $n \in \omega$
so that $\diam(F^n) \le 2^{-n}$ and each
$F^{n+1} \subsetneqq F^n$.  Let $\{p_\alpha\} = \bigcap_n F^n$
and let $q_\alpha^n$ be any point in $F^{n+1} \backslash F^n$.
Make sure that $F^0 \subseteq Q^0 = P_\alpha$ whenever 
$P_\alpha$ is uncountable, so that $p_\alpha \in P_\alpha$ is as
required by (9).  Also make sure that for every $n$,
either $F^n \subseteq Q^n$ or $F^n \cap Q^n = \emptyset$,
so that (8) will hold.
\end{proofof}

\section{Peak Sets}
\label{sec-peak}
Fix $\alpha < \omega_1$.
The function $h_\alpha$ occurring in the proof of Theorem \ref{thm-weird-HL}
is easy to construct because $X_\alpha$ is a compact metric space.
Note that there are also uniformly bounded $g_{\alpha,n} \in C(X_\alpha)$
(for $n \in \omega$)
with $g_{\alpha,n}(x) \to h_\alpha(x)$ whenever $x \ne p_\alpha$.
In the proof of Theorem \ref{thm-main}, we shall furthermore require that
each $g_{\alpha,n} \in \AA_\alpha$, where $\AA_\alpha \leac C(X_\alpha)$.
This is not always possible.  For example, if $X_\alpha = \cdisc$
and $\AA_\alpha = \DDDD$, the disc algebra, then we could not
find such $g_{\alpha,n}$ and $h_\alpha$ unless $p_\alpha \in \TTT$,
since $h_\alpha$ is required to be discontinuous at $p_\alpha$.
For $\alpha= 1$, we shall avoid this problem by defining
$X_1$ to be $\TTT$; then
a suitable $h_1$ can be concocted using standard facts about
$H^\infty$ (see \cite{GA, HO, RR, RUD3}).
To obtain suitable $h_\alpha$ on $X_\alpha$ for $\alpha > 1$,
we shall require that all points of $X_\alpha$ be peak points;
the following is easily seen to be equivalent to the usual definition
(see, e.g., \cite{GAM}):

\begin{definition}
\label{def-peak}
Assume that $X$ is compact, $\AA \lea C(X)$, and $H$ is a closed subset 
of $X$.  Then $H$ is a \emph{peak set}
\textup(with respect to  $\AA$\textup) iff
there is an $f \in \AA$ such that
\begin{itemizz}
\item[1.] $f(x) = 0$ for all $x \in H$.
\item[2.] $\Re(f(x)) > 0$ for all $x \notin H$.
\end{itemizz}
$\PS_\AA (X)$ is the set of all $H \subseteq X$
which are peak sets with respect to $\AA$.
$p \in X$ is a \emph{peak point} iff $\{p\}$ is a peak set.
\end{definition}

Every peak set is a closed $G_\delta$ set, but not conversely.
For example, if $H$ is clopen and $\AA \leac C(X)$,
then by Runge's Theorem, $H$ is a peak set iff
$\cchi_H \in \AA$. Also, for the disc algebra,
$p \in \cdisc$ is a peak point iff $|p| = 1$.

Our primary interest here is in the peak \emph{points}.
However, we mention peak \emph{sets} because these will be used to
prove that $\PS_\AA (X)$ contains singletons
by applying the following well-known fact:

\begin{lemma}
\label{lemma-peak-int}
If $\AA \leac C(X)$, then $\PS_\AA (X)$ is closed under
countable intersections and finite unions.
\end{lemma}
\begin{proof}
For intersections, fix $H_n \in \PS_\AA (X)$ for
$n \in \omega$, and let $H = \bigcap_n H_n$.
Let $f_n$ satisfy (1)(2) of Definition \ref{def-peak}
for $H_n$, and assume that $\|f_n\| \le 2^{-n}$.
Let $f = \sum_n f_n$.  Then $f \in \AA$ because $\AA$ is closed,
and $f$ satisfies (1)(2) for $H$.

For unions, let $H = H_0 \cup H_1$, and let 
$f_0, f_1$ satisfy (1)(2) of Definition \ref{def-peak}
for $H_0, H_1$ respectively.
Define $f(x) = \sqrt{ f_0(x)}\sqrt{ f_1(x)}$.
Again, $f \in \AA$ because $\AA$ is closed,
since $\sqrt z $ can be uniformly approximated by polynomials
on any compact subset of $\{z\in\CCC : \Re(z) \ge 0\}$, and
$f$ satisfies (1)(2) for $H$.
\end{proof}

\begin{lemma}
\label{lemma-peak-limit}
Assume that $X$ is compact,
$\AA \leac C(X)$, and $p\in X$ is a peak point.
Let $\langle q_n : n \in \omega \rangle$ be a sequence
of points in $X \backslash \{p\}$ converging to $p$.  Then there
are functions $h$ and $g_n$ for $n\in \omega$ such that:
\begin{itemizz}
\item[1.] Each $g_n \in \AA$.
\item[2.] Each $\|g_n\| \le 1$.
\item[3.]  $h \in C(X \backslash \{p\}, \cdisc)$.
\item[4.]  On $X \backslash \{p\}$, the $g_n$ converge
to $h$  uniformly on compact sets.
\item[5.] $|h(x)| \to 1$ as $x \to p$ in $X \backslash \{p\}$.
\item[6.] Every $w \in \TTT$ is a limit point of
the sequence $\langle h(q_n) : n \in \omega \rangle$.
\end{itemizz}
\end{lemma}

\begin{proof}
Let $f_0$ be the function given by Definition \ref{def-peak}.
We plan to obtain $h$ by composing $f_0$ with a suitable Blaschke product.
The notation will be easier if we define the product
in the upper halfplane; see, e.g., \cite{GA},\S II.2.
Let 
\[
V = \{x + iy \in \CCC \ :\  0 < -x < y    \} \ \ .
\]
If $f(z) =   e^{5\pi i / 8} \cdot \sqrt[4]{f_0(z)}$,
then $f \in \AA$, $f(p) = 0$, and $f(x) \in V$ for all $x \ne p$.
When $\Im(\alpha) > 0$, let
\[
B_\alpha(z) = \frac{z - \alpha}{z - \overline\alpha} \ \ .
\]
Then $|B_\alpha(z)|$ is $1$ on the real axis and less than $1$
in the upper halfplane.  Let $z_\ell = f(q_\ell) \in V$; then
$z_\ell \to 0$.  We shall choose
$\alpha_n$ in the upper halfplane and form the Blaschke products:
\[
B^{(n)} (z) = \prod_{m<n} B_{\alpha_m}(z) \qquad\quad
B(z) = \prod_{n\in\omega} B_{\alpha_n}(z) 
\]
They will satisfy:
\begin{itemizz}
\item[a.] $B^{(n)}(z) \to B(z)$ uniformly on compact subsets of $V$.
\item[b.] $|B(z)| \to 1$ as $z \to 0$ in $V$.
\item[c.] Every point in $\TTT$ is a limit point of the sequence
$\langle B(z_\ell) : \ell \in \omega \rangle$.
\end{itemizz}
Assuming that this can be done, the lemma is satisfied by
letting $g_n = B^{(n)} \circ f$ and $h = B \circ f$.
$g_n \in \AA$ because 
each $B^{(n)}$ is holomorphic in a convex neighborhood of 
$f(X)$, and hence can be uniformly approximated
on $f(X)$ by polynomials.

To obtain (a)(b)(c), we choose the $\alpha_n$,
along with a subsequence, 
$\langle z_{\ell_n}: n \in \omega \rangle$, of
$\langle z_\ell: \ell \in \omega \rangle$,
to satisfy:
\begin{itemizz}
\item[d.] $\alpha_n = \xi_n + i \eta_n$ and $0 < \xi_n = (n+1) \eta_n$.
\item[e.] $z_{\ell_n}  = x_n + iy_n$ and $\eta_n = y_n$.
\item[f.] $n > m \Rightarrow \xi_n \le  2^{-n}\eta_m$.
\item[g.] $n > m \Rightarrow | 1 -
\arg(B_{\alpha_m}(z_{\ell_n})) / \arg(B_{\alpha_m}(0)) | \le 2^{-n}$.
\end{itemizz}
So, $\alpha_n$ is to the right of $V$ and 
$ \xi_0 \ge \eta_0 \ge 
\xi_1 \ge \eta_1 \ge 
\xi_2 \ge \eta_2 \ge \cdots $.
The $\alpha_n$ and $z_{\ell_n}$ can easily be chosen by induction
to satisfy (d)(e)(f)(g),
using $z_\ell \to 0$ and the continuity of $B_{\alpha_m}$ at $0$.
We now verify (a)(b)(c).
Observe that $\eta_n / \xi_n \to 0$ but
$\sum_n \eta_n / \xi_n = \infty$; this will allow us to prove (b)
without having $\lim_{z \to 0} B(z)$ exist, which would contradict (c).

For (a), note that if $\alpha = \xi + i\eta$ and $z = x + iy$ then
\[
B_\alpha(z) = \frac{x - \xi + iy  - i\eta } {x- \xi  + iy  + i\eta } = 
1 - \frac{ 2 i\eta } {x- \xi  + iy  + i\eta } \ \ . 
\]
Then, as usual with Blaschke products, (a) follows from
$\sum_n \eta_n < \infty$, which in turn follows from (d)(f).

For (b),
we need to estimate $|B_\alpha(z)|$, where 
$\alpha = \xi + i \eta$, 
$z = x + i y \in V$, and $0 < \eta \le \xi$.
Now
\[
|B_\alpha(z)|^2 = 
\frac{ (\xi - x)^2 + (y - \eta)^2 }
{(\xi - x )^2 + (y  + \eta)^2 } =
1 - \frac{4y\eta } {(\xi - x )^2 + (y  + \eta)^2 } 
\ \ . 
\]
Clearly,
\[
1 > |B_\alpha(z)|^2 \ge
1 - \frac{4\eta }{y} \qquad 
1 > |B_\alpha(z)|^2 \ge
1 - \frac{4y }{\eta}  \ \ ,  
\tag{$\ast$}
\]
and these are useful when $\eta \ll y$ or $y \ll \eta$.  But also note that:
\[
1 > |B_\alpha(z)|^2 \ge 1 -  \frac{4 \eta}{\xi}  \ \ .
\tag{\dag}
\]
To prove this:  If $y \ge \xi$ then (\dag) follows from $(\ast)$.
If $y \le \xi$ then, since $x \le 0$, 
$|B_\alpha(z)|^2 \ge 1 - (4y\eta)/(\xi^2) \ge 1 - (4\eta)/\xi$.

To prove (b), fix $z = x + i y \in V$ with $y \le \eta_1$.
Next fix $n\ge 1$ such that $\eta_{n+1} \le y \le \eta_n$.
We show that $|B(z)| = 1 - o(1)$ as $n \nearrow \infty$ by
estimating each $|B_{\alpha_m}(z)|^2$.
Applying (\dag) and (d), we get
$|B_{\alpha_n}(z)|^2 \ge 1 - 4 \eta_n / \xi_n = 1 - 4/(n+1)$,
and likewise $|B_{\alpha_{n+1}}(z)|^2 \ge 1 - 4/(n+2)$.
For $m < n$ use $(\ast)$ and (d)(f) to get
\[
1 > |B_{\alpha_m}(z)|^2 \ge 1 - \frac{4y }{\eta_m} 
\ge 1 - \frac{4\eta_n }{\eta_m} 
= 1 - \frac{4\xi_n }{(n+1)\eta_m} 
\ge 1 - \frac{4\cdot 2^{-n}}{n + 1}  
\ \ ,  
\]
so $\prod_{m < n} |B_{\alpha_m}(z)|^2 \ge 1 - 4\cdot 2^{-n}$.
For $m > n+1$ use $(\ast)$ and (d)(f) to get
\[
1 > |B_{\alpha_m}(z)|^{2} \ge
1 - \frac{4\eta_m }{y} \ge
1 - \frac{4\eta_m }{\eta_{n+1}} \ge
1 - 4 \cdot 2^{-m}
\ \ ,
\]
so $\prod_{m > n+1} |B_{\alpha_m}(z)|^{2} \ge 1 - 2^{-n+1}$.
Putting these estimates together, we get $|B(z)| = 1 - o(1)$.

To verify (c), note that
(b) implies that (c) is equivalent to the assertion that
$\{\arg(B(z_{\ell_n})) : n \in \omega\}$ is dense in $\TTT$.
We estimate $\arg(B(z))$, using:
\[
\arg(B_\alpha(z)) =
\arctan\frac{\eta - y}{\xi - x} +
\arctan\frac{\eta + y}{\xi - x} =
\arctan\frac{ 2  \eta (\xi - x )}{(\xi -  x )^2  + y^2 - \eta^2 }
\ \ . 
\]
We are using
$\arctan(u) + \arctan(v) = \arctan( (u+v) / (1 - uv))$;
this applies here because all three of 
$\arg(\alpha - z)$,
$\arg(\overline\alpha - z)$, and $\arg(B_\alpha(z))$ are in 
the range $(-\pi/2, \pi/2)$.
Let $\theta^m_n =   \arg(B_{\alpha_m}(z_{\ell_n}))$.
Then $\arg(B(z_{\ell_n})) \equiv \sum_m \theta^m_n \bmod 2\pi$.
Define:
\[
\sigma^m := 
  \arg(B_{\alpha_m}(0)) =
2 \arctan \frac{\eta_m }{\xi_m} =
 2 \arctan \frac{1 }{(m+1)}
\ \ ,
\]
and observe that we have:
\begin{itemizz}
\item[1.] $0 < \sigma^m \to 0$ and $\sum_m \sigma^m = \infty$.
\item[2.] $\theta^n_n \to 0$.
\item[3.] $\sum_{m > n} |\theta^m_n|  \le 2^{-n+2}$.
\item[4.] $m < n \Rightarrow \sigma^m(1 - 2^{-n}) \le
\theta^m_n  \le \sigma^m(1 + 2^{-n})$.
\end{itemizz}
(1) holds because $\sigma^m \approx 2/m$.
(2) follows from  (e) and $x_n \le 0$, which yields
\[
\theta^n_n  =   \arctan\frac{2\eta_n}{\xi_n - x_n} =
  \arctan\frac{2\eta_n}{(n+1)  \eta_n - x_n} \le
  \arctan\frac{2\eta_n}{(n+1)  \eta_n } \to 0 \ \ .
\]
For (3), use  (d)(f) to get, for $m > n$:
\[
|\theta^m_n| =  
\arctan\frac{ 2  \eta_m |\xi_m - x_n |}
{(\xi_m -  x_n )^2  + \eta_n^2 - \eta_m^2 } \le
 \arctan\frac{ 4  \eta_m \eta_n} {\eta_n^2  } \le
\frac{ 4  \xi_m } {(m + 1)\eta_n  } \le 4 \cdot 2^{-m} \ \ ,
\]
so that
$\sum_{m > n} |\theta^m_n|  \le 4  \cdot2^{-n}$.
\;(4) is immediate from (g).

Finally, (1) implies that the values
$ \sum_{m < n} \sigma_m \bmod 2\pi$ (for $n\in\omega$) are dense in $\TTT$,
and (2)(3)(4) imply that as $n \to \infty$, these
values get close to $\arg(B(z_{\ell_n}))$.
\end{proof}

We remark that there are well-known 
interpolation theorems of Pick, Nevanlinna, Carleson, and others
(see \cite{GA, HO, RR}) which involve constructing Blaschke products 
to have given values on a given sequence of points.  However,
because of our requirement (b) in the above proof, we do not see
how to obtain our Blaschke product simply by
quoting one of these theorems.

\section{Subspaces of Polydiscs}
\label{sec-poly}
We now return to the construction of \S \ref{sec-construct},
and show how to modify the space so that it also fails the CSWP.
To get a function algebra witnessing this failure, it is easier
to construct the space in $\cdisc\,^{\omega_1}$ rather than $[0,1]^{\omega_1}$,
so we start by replacing $[0,1]$ with $\cdisc$ 
in the requirements of \S  \ref{sec-construct}.

We shall get
$X  = X_{\omega_1} \subseteq \cdisc\,^{\omega_1}$, with
$X_\alpha = \pi^{\omega_1}_\alpha(X) \subseteq \cdisc\,^{\alpha}$.
Let $\RF^-$ denote the requirements consisting of
conditions (1)--(5) and (8)--(9) of 
\S \ref{sec-construct} plus:
\begin{itemizz}
\rm
\item[$\widetilde 0$.] $X_1 = \TTT$.
\item[$\widetilde 6$.] $p_\alpha \in X_\alpha$ and
$h_\alpha \in C(X_\alpha \backslash \{p_\alpha\},\; \cdisc)$
and $X_{\alpha+1} = \overline{h_\alpha}$.
\item[$\widetilde 7$.] $q_\alpha^n \in  X_\alpha \backslash\{p_\alpha\}$, and
$\langle q_\alpha^n: n \in \omega \rangle \to p_\alpha$,
and all points of $\TTT$ are limit points of
$\langle h_\alpha(q_\alpha^n): n \in \omega \rangle$, and
$|h_\alpha(x)| \to 1$ as $x \to p_\alpha$ in $X_\alpha \backslash \{p_\alpha\}$,
and $\{p_\alpha\} \times \TTT \in \PP_{\alpha+1}$.
\end{itemizz}

Note that we have the slice $\{z \in \cdisc : (p_\alpha, z) \in X_{\alpha+1}\}$
equal to $\TTT$, not $\cdisc$, as one might expect.
This will enable us to prove that all points in each $X_\beta$
are peak points; see Lemma \ref{lemma-peak}.
Since $\TTT$ is connected, the argument is essentially unchanged,
and we get a weird HL space as before, using $\diamondsuit$.

Along with the $X_\alpha$, we need a function
algebra on $X_\alpha$ refuting the CSWP.
We use the obvious analog of the disc algebra:
\begin{definition}
$\PPPP_\alpha \lea C(\cdisc\,^\alpha)$ is the algebra
generated by the projections
$\{ \varphi_\xi : \xi < \alpha\}$
\textup(see Definition \ref{def-proj}\,\textup), and
$\DDDD_\alpha \leac C(\cdisc\,^\alpha)$ is the uniform
closure of $\PPPP_\alpha $.
Let $\AA_\alpha$ be the uniform closure of $\PPPP_{\alpha}\res X_\alpha %
= \{f \res X_\alpha: f \in \PPPP_\alpha\}$
\end{definition}

For finite $\alpha$, $\PPPP_\alpha$ is the algebra
of polynomials in $\alpha$ complex variables on the 
polydisc $\cdisc\,^\alpha$, and 
$\DDDD_\alpha$ is the algebra of continuous functions
which are holomorphic in the interior of the polydisc.
For all $\alpha > 0$,
$\DDDD_\alpha \ne C(\cdisc\,^\alpha)$.
In constructing the $X_\alpha$, we also make sure that
$\AA_\alpha  \ne C(X_\alpha)$.
To do this, we choose all $h_\alpha$ in $H^\infty$.
More precisely:

\begin{definition}
Let $\lambda = \lambda_1$ be the Haar probability measure on
$X_1 = \TTT$.  For $1 \le \alpha < \omega_1$, let 
$\lambda_\alpha$ be the unique Borel probability measure
on $X_\alpha$ such that $\lambda_1$ is the induced measure
$\lambda_\alpha\,(\pi^\alpha_1)\iv$.
For $1 \le \alpha \le \beta < \omega_1$, define the map
$(\pi^\beta_\alpha)^* : L^\infty(X_\alpha, \lambda_\alpha) \to 
L^\infty(X_\beta, \lambda_\beta) $ by
$(\pi^\beta_\alpha)^*([f]) = [f \circ \pi^\beta_\alpha]$,
where
$[g] \in L^\infty$ denotes the equivalence class of $g$.
\end{definition}
Note that each $\lambda_\alpha$ is unique because all points in $X_1$
outside the countable $\{\pi^\xi_1(p_\xi) : 1 \le \xi < \alpha\}$
have a unique preimage under $\pi^\alpha_1$.  Likewise,
$(\pi^\beta_\alpha)^*$ is a Banach algebra isomorphism.

Let $\RF$ consist of the requirements of $\RF^-$,
along with this requirement on the $h_\alpha$:

\begin{itemizz}
\item[$\widetilde{10}$.] For $1 \le \alpha < \omega_1$,
$[h_\alpha] \in (\pi^\alpha_1)^*(H^\infty(\TTT))$.
\end{itemizz}

This makes $X$ fail the CSWP.  
Requirement $(\widetilde{10})$ is used explicitly
in the proof of the next lemma. 
Lemma \ref{lemma-not-dense} follows, and
produces a continuous function
not in $\AA_{\omega_1}$.

\begin{lemma}
\label{lemma-P-H}
Fix $\beta$ with $1 \le \beta < \omega_1$.
Suppose requirement $(\widetilde{10})$ holds for all $\alpha < \beta$.
Then $[k] \in (\pi^\beta_1)^*(H^\infty(\TTT))$
for each $k \in \AA_\beta$. 
\end{lemma}
\begin{proof}
Since $\PPPP_{\beta}\res X_\beta$ is generated by
$\{\varphi_\alpha : \alpha < \beta\}$, it suffices to
prove that each $[\varphi_\alpha] \in (\pi^\beta_1)^*(H^\infty(\TTT))$.
Now, $[\varphi_0] = (\pi^\beta_1)^*([I])$, where $I(z) = z$.
For $1 \le \alpha < \beta$,
$[\varphi_\alpha] = (\pi^\beta_{\alpha})^*([h_\alpha]) %
= [h_\alpha \circ \pi^\beta_\alpha]$. 
By $(\widetilde{10})$ for $\alpha < \beta$,
$[h_\alpha] = [h \circ \pi_1^\alpha]$ for some $h \in H^\infty(\TTT)$.
So $[\varphi_\alpha] = [h \circ (\pi_1^\alpha \circ \pi^\beta_\alpha)] %
\in (\pi^\beta_1)^*(H^\infty(\TTT))$.
\end{proof}

\begin{lemma}
\label{lemma-not-dense}
Suppose requirements $\RF$ hold.
Then $\AA_{\omega_1}  \ne C(X)$.
\end{lemma}
\begin{proof}
Let $\overline I \in C(X_1)$ denote the usual complex conjugation
given by $\overline I (z) = \overline z$.
Then $\overline I \circ \pi^{\omega_1}_1$
(i.e., $\vec z \mapsto \overline{z_0}$) is not in $\AA_{\omega_1}$.
To see this: it suffices to show that
$\overline I \circ \pi^{\beta}_1 \notin \AA_\beta$
for all $\beta < \omega_1$.
Since $\overline I \notin H^\infty(\TTT)$, 
$(\pi^{\beta}_1)^*([\overline I]) \notin (\pi^{\beta}_1)^*(H^\infty(\TTT))$
for all $\beta < \omega_1$.
So the result follows from Lemma \ref{lemma-P-H}. 
\end{proof}

\begin{lemma}
\label{lemma-peak}
Fix $\beta$ with $1 \le \beta < \omega_1$.
Suppose requirement $(\widetilde{10})$ holds for all $\alpha < \beta$.
Then each $y \in  X_\beta$  is a peak point with respect to  $\AA_\beta$.
\end{lemma}
\begin{proof}
We induct on $\beta$.
For $\beta = 1$, this is clear, since $X_1 = \TTT$.

If $\beta$ is a limit, then $\{y\} = 
\bigcap_{\alpha < \beta} (\pi^\beta_\alpha)\iv  \{(\pi^\beta_\alpha)(y)\} $.
Applying the lemma inductively, each 
$(\pi^\beta_\alpha)(y)$ is a peak point in $X_\alpha$
with respect to $\AA_\alpha$, which implies that
each $ (\pi^\beta_\alpha)\iv  \{(\pi^\beta_\alpha)(y)\} $
is a peak set in $X_\beta$ with respect to $\AA_\beta$.
The result now follows using Lemma \ref{lemma-peak-int}.

Now, say $\beta = \alpha+1$, let $v = \pi^\beta_\alpha(y)$
and let 
$H =  (\pi^\beta_\alpha)\iv  \{v\} $,
which, as above, is a peak set in $X_\beta$.
If $v \ne p_\alpha$, then $H = \{y\}$.
If $v = p_\alpha$, then $y \in H = \{v\} \times \TTT$ 
(using condition $(\widetilde 7))$.
If $y = (v, e^{i\theta})$, then 
$K = \{x \in X_\beta:  \varphi_\alpha(x) = e^{i\theta} \}$
is also a peak set, and $\{y\} = H \cap K$.
\end{proof}

\begin{proofof}{Theorem \ref{thm-main}}
We need to show inductively that requirements $\RF$ can indeed be met.
Suppose that we have constructed $X_\beta$ so that
they hold for all $\alpha < \beta$.
Get $p_\beta\in X_\beta$ and $\langle q_\beta^n : n \in \omega \rangle$
converging to $p_\beta$ 
as in the proof of Theorem \ref{thm-weird-HL}.
By Lemma \ref{lemma-peak}, $p_\beta$ is a peak point.
Now get $h \in C(X_\beta \backslash \{p_\beta\}, \cdisc)$ and
$g_n \in \AA_\beta$ as in Lemma \ref{lemma-peak-limit}.
Then each $[g_n] \in (\pi^\beta_1)^*(H^\infty(\TTT))$
by Lemma \ref{lemma-P-H}, so
$[h] \in (\pi^\beta_1)^*(H^\infty(\TTT))$ since 
$g_n \to h$ on $X_\beta \backslash \{p_\beta\}$.
Thus, taking $h_\beta = h$ satisfies
$(\widetilde{10})$ for $\beta$.
Lemma \ref{lemma-peak-limit} also guarantees that this choice
of $h_\beta$ will satisfy the rest of $(\widetilde{7})$.
The remaining requirements are satisfied as for 
Theorem \ref{thm-weird-HL}.
\end{proofof}

\section{Some Forcing Orders}
\label{sec-forcing}

\begin{definition}
Order $2^{<\omega_1}$ by: $p \le q$ iff $p \supseteq q$.
Let $\one = \emptyset$, the empty sequence.
\end{definition}

So, $2^{<\omega_1}$ is a tree, with the root $\one$ at the top.
Viewed as a forcing order, it is equivalent to countable partial
functions from $\omega_1$ to $2$.  This forcing is countably closed,
and thus preserves all witnesses to $\diamondsuit$, and thus preserves
the weird space constructed in the proof of Theorem
\ref{thm-weird-HL}.
To kill such spaces, we shall force with subtrees of 
$2^{<\omega_1}$ which satisfy a weakening of countable closure.

\begin{definition}
A \emph{Cantor tree of sequences} is a subset
$\{p_s : s \in 2^{<\omega}\} \subseteq 2^{<\omega_1}$ such
that each $p_{s\cat\mu} < p_s$ for $\mu = 0,1$, 
and $p_{s\cat 0}\perp p_{s\cat 1}$.
\end{definition}
That is, $p_{s\cat 0}$  and $p_{s\cat 1}$ are incompatible extensions of $p_s$.

\begin{definition}
$\PPP \subseteq 2^{<\omega_1}$ has the \emph{Cantor tree property}
iff:
\begin{itemizz}
\item[1.] $\one \in \PPP$ and
$\PPP $ is a subtree:  $q \ge p \in \PPP \to q \in \PPP$.
\item[2.] If $p \in \PPP$ then $p\cat 0 , p \cat 1 \in \PPP$.
\item[3.] Whenever $\{p_s : s \in 2^{<\omega}\} \subseteq \PPP$
is a Cantor tree of sequences, there is at least one
$f \in 2^\omega$ such that $\bigcup\{p_{f\res n} : n \in \omega\} \in \PPP$.
\end{itemizz}
\end{definition}
Of course, then by (3) there must be uncountably many such $f$;
in fact the set of $f$ satisfying (3) must meet every perfect subset
of the Cantor set $2^\omega$, since otherwise we could find a subtree
of the given Cantor tree which contradicts the Cantor tree property.
It is also easily seen by induction that $\PPP$ is a normal subtree; 
i.e.:

\begin{lemma}
\label{lemma-normal}
If $\PPP \subseteq 2^{<\omega_1}$ has the Cantor tree property, then
whenever $p \in \PPP$ and $\dom(p) < \alpha < \omega_1$, there is 
a $q \in \PPP \cap 2^\alpha$ such that $q < p$.
\end{lemma}

If $\PPP$ has the Cantor tree property,
then it is proper and forcing with it adds no $\omega$-sequences.
Such orders are called \emph{totally proper};
see Eisworth and Roitman \cite{ER}, which gives a number of equivalents,
which we use in:

\begin{lemma}
\label{lemma-tot-prop}
If  $\PPP \subseteq 2^{<\omega_1}$ has the Cantor tree property, then
$\PPP$ is totally proper.
\end{lemma}
\begin{proof}
Fix a suitably larger regular cardinal, and let $M \prec H(\theta)$
be countable and fix $p \in \PPP \cap M$.  Following \cite{ER},
it is sufficient to find a $q \le p$ such that whenever $A \subseteq \PPP$
is a maximal antichain and $A \in M$, there is an $r \in A\cap M$ with
$q \le r$.

To get $q$, let $\{A_n : n \in \omega\}$ list all the maximal
antichains which are in $M$, and build a Cantor tree
$\{p_s : s \in 2^{<\omega}\} \subseteq \PPP \cap M$ such that
$p_\emptyset \le p$ and $p_s$ extends an element
of $A_n \cap M$ for each $s \in 2^n$.  Then choose
$f \in 2^\omega$ such that
$q := \bigcup\{p_{f\res n} : n \in \omega\} \in \PPP$.
\end{proof}

Thus, assuming PFA, this $\PPP$ will have an uncountable
chain.  By Lemma \ref{lemma-getP}, a weird space will yield such a $\PPP$,
and hence cannot be HL under PFA.

\begin{lemma}
\label{lemma-con-split}
If $X$ is compact, connected, and infinite, and  $U \subseteq X$
is a nonempty open set, then there is a closed $K \subseteq U$
such that $K$ is connected and infinite.
\end{lemma}
\begin{proof}
Let $V$ be open and nonempty with $\overline V \subseteq U$,
fix $p \in V$, and let $K = \comp(p, \overline V)$.
If $K = \{p\}$, then there is an $H$ which is relatively 
clopen in $\overline V$ such that $H \subseteq V$.  But then $H$ would
be clopen in $X$, a contradiction.
\end{proof}

\begin{lemma}
\label{lemma-getP}
Assume that $X$ is compact, HL, and not totally disconnected,
and assume that $X$ has no 
subspace homeomorphic to the Cantor set $2^\omega$.  Then
there exists a $\PPP$ with the Cantor tree property 
which has no uncountable chains.
\end{lemma}
\begin{proof}
Along with $\PPP$, we shall choose sets $H_p$ for $p \in \PPP$
with the following properties:
\begin{itemizz}
\item[1.] $H_p$ is an infinite closed connected subset of $X$.
\item[2.] If $p \in \PPP$, then $p\cat0, p\cat1 \in \PPP$
and $H_{p\cat0}, H_{p\cat0}$ are disjoint subsets of $H_p$.
\item[3.] If $p \in 2^\gamma$, where $\gamma$ is a countable limit ordinal
and $p\res\alpha \in \PPP$ for all $\alpha < \gamma$, then
$H_p = \bigcap\{H_{p\res\alpha} : \alpha < \gamma\}$, and 
$p \in \PPP$ iff $H_p$ is infinite.
\end{itemizz}
$H_\one$ can be chosen because $X$ is not totally disconnected.
Given $p\in \PPP$, we can choose $H_{p\cat0}, H_{p\cat1}$ by applying
Lemma \ref{lemma-con-split} to $H_p$.
To verify the Cantor tree property,
let $\{p_s : s \in 2^{<\omega}\} \subseteq \PPP$
be a Cantor tree of sequences.  For $f \in 2^\omega$,
let $p_f = \bigcup\{p_{f\res n} : n \in \omega\}$.
If none of these $p_f$ are in $\PPP$, then each
$H_{p_f}$ would be a singleton, $\{x_f\}$.
But then  $\{x_f : f\in 2^\omega\}$ is homeomorphic to the Cantor set.
\end{proof}

\begin{proofof}{Theorem \ref{thm-PFA-HL}}
By Lemmas \ref{lemma-tot-prop} and \ref{lemma-getP}.
\end{proofof}

\end{document}